\documentclass{amsart}
\usepackage{url}

\def\fr#1/#2{\frac{#1}{#2}}				% Display fraction
\def\tfr#1/#2{{#1}/{#2}}				% Text fraction

\author{Jacques G\'elinas}
\address{Ottawa, Canada}

\email{jacquesg00@hotmail.com}

\thanks{This work was done in December 2016, while the author was a retired mathematician}

\begin{document}

%=================================================================================================

\keywords{Stirling formula, Bernoulli numbers, Wallis product, De Moivre table of logarithms}
\subjclass{Primary 41A60, Secondary 11B68}

\title{Original proofs of Stirling's series for log(n!)}

%\maketitle
%\vspace{-2\baselineskip}
%\begin{center} \footnotesize 2017-01-20 \end{center}

\begin{abstract}
This is a transcription into modern notations of the derivation by Stirling and De Moivre
of an asymptotic series for $\log(n!)$, usually called Stirling's series. The previous
discovery by Wallis of an infinite product for $\pi$, and later results on the divergence
of the series are also presented. We conclude that James Stirling has priority over
Abraham de Moivre for Stirling's formula and Stirling's series.
\end{abstract}

\maketitle

%=================================================================================================

\section{Origin of Stirling's series}

In order to evaluate the central binomial coefficient $\binom{2n}{n} \approx \fr {2^{2n}}/{\sqrt{\pi n}}$
for his research in probability theory \cite{DeMoivre:1718}, Abraham De Moivre needed a formula for $\log n!$ to replace
the table of 14 figure logarithms of factorials from 10! to 900! by differences of 10 printed in his new book
``Miscellanea Analytica" \cite{DeMoivre:1730}.
A ``few days" after publication in 1730, James Stirling sent him a letter pointing out numerical errors in that
table past the fifth decimal making it ``unsuitable for research", and stated without proof a series for 
$\log{n!}$, involving the variable $z=n+\tfr 1/2$, the constant $\log\sqrt{2\pi}$, and certain rational 
coefficients defined by a triangular system of linear equations~:
\begin{equation} \label{Sform}
\log n! = z\log z - z + \log\sqrt{2\pi} - \fr 1/{2 \times 12 z} + \fr 7/{8 \times 360 z^3}
					- \fr 31/{32 \times 1260 z^5} + \ldots
\end{equation}
De Moivre saw immediately as factors in these coefficients the Bernoulli numbers $B_{2k}$ that he knew 
$(12=1\cdot2\cdot6,\,360=3\cdot4\cdot30,\,1260=5\cdot6\cdot42)$ and proceeded to derive independently a 
closely related series involving more simply the variable $n$ itself~:
\begin{equation} \label{Dform}
\log n! = \left(n+\fr 1/2\right)\log n - n + 1 + 
			\sum_{k=1}^\infty \fr B_{2k}/{2k(2k-1)}\left(\fr 1/{n^{2k-1}} - 1 \right).
\end{equation}
Stirling gave the proof of (\ref{Sform}) in his book ``Methodus Differentialis" that same year \cite{Stirling:1730}.
De Moivre later published the very different proof of (\ref{Dform})
through a 22 page ``Supplementum" \cite{DeMoivre:1731} to ``Miscellanea Analytica", giving due credit to Stirling
by quoting explicitly (\ref{Sform}), and even appending it in 1756 to his book
``Doctrine of chances" \cite{DeMoivre:1718}, along with a corrected table of logarithms (see below).

Both series are usually called Stirling's series, but the original one by Stirling
has an explicit constant and can be used, while the simpler modified one by De Moivre involves a constant 
expression whose value must still be determined by other means, for example from Wallis infinite product.

In summary, the priority for the discovery of the series and the determination of the constant
$\log \sqrt{2\pi}$ was given to Stirling by De Moivre himself, while the modifications in its expression 
and in its proof from the recognition of Bernoulli numbers must be attributed to De Moivre.

Unfortunately, about half of the known copies of the ``Miscellanea Analytica"
lack the ``Supplementum" (translated in \cite{LanierTrotoux:1996}), and some publications 
or Wikipedia web pages unaware of its contents attribute incorrectly to De Moivre the discovery
of the series, leaving only the determination of the constant to Stirling's credit
\cite{Namias:1986,Conrad:2014}.
However, this does not diminish the importance of De Moivre in the early history of 
statistics, which is now well documented and accepted \cite{Pearson:1924}.

\section{Derivation of the simpler series by De Moivre}

First, $\log n!$ is cleverly expressed as a linear combination of $n$ simpler logarithms.
\begin{align*}
	\log n! &= \log n + \log \left( n^{n-1}\right) + \log (n-1)! - \log \left( n^{n-1}\right)
\\
	&= \log n + (n-1)\log n + \log \prod_{k=1}^{n-1} \left( \fr {n-k}/n \right)
\\
	&= n\log n + \sum_{k=1}^{n-1} \log \left( 1 - \fr k/n \right).
\end{align*}

Next, each logarithm in the finite sum of $n-1$ logarithms is replaced by an infinite Taylor series.
$$
	 - \log \left( 1 - \fr k/n \right) =  \sum_{p=1}^\infty \fr 1/p \left( \fr k/n \right)^p
		=  \sum_{p=1}^\infty \fr 1/ p \fr 1/{n^p} k^p.
 $$

Summing this on $k$ produces a double summation involving, after permutation, sums of powers of the first integers
for which James Bernoulli had published a nice explicit formula earlier in 1713 involving the factors
$B_k=1, \fr 1/2, \fr 1/6, 0, -\fr 1/{30}, 0, \fr 1/{42}, 0, -\fr 1/{30}, 0, \fr 5/{66}, \ldots$
that De Moivre called ``Bernoulli numbers\footnote{The sum up to $n-1$ has $B_1=+\fr 1/2$ in terms of $n-1$,
instead of the usual $B_1=-\fr 1/2$ in terms of $n$.}". 
After extracting the first two terms, De Moivre inverted the order of summation yet again and derived
as examples the first four terms of the infinite series. We also set with him $x=\tfr (n-1)/n=1-\tfr 1/n$
but, after inserting one term, modify one index in order to proceed along diagonals and thus obtain the general
expression.
\begin{align*}
	\sum_{p=1}^\infty \fr 1/p \fr 1/{n^p}  \sum_{k=1}^{n-1} k^p
	 &= \sum_{p=1}^\infty  \fr 1/p \fr 1/{n^p} \left\{\fr 1/{p+1}\sum_{k=0}^p\binom{p+1}{k} B_k (n-1)^{p+1-k}\right\}
\\
	&= \sum_{p=1}^\infty \left( \fr {n x^{p+1}}/{p(p+1)} + \fr x^{p}/{2p} \right)
	+ \sum_{p=2}^\infty \sum_{k=2}^p\fr 1/{p(p+1)}\binom{p+1}{k} \fr {B_k}/{n^{k-1}} x^{p+1-k}
\\
	&= \int_0^x\sum_{p=1}^\infty \left( n\fr {t^p}/p + \fr t^{p-1}/2 \right) \, dt
	+ \sum_{k=2}^\infty \fr {B_k}/{n^{k-1}} 
			\left[ \sum_{p=k}^\infty  \fr 1/{p(p+1)} \binom{p+1}{p-k+1} x^{p-k+1} \right]
\\
	&= \int_1^{1-x}\left( n\log t + \fr 1/{2t} \right)\,dt
	+ \sum_{k=2}^\infty \fr {B_k}/{n^{k-1}} \left[ \sum_{r=0}^\infty  \fr \binom{r+k}{k}/{(r+k)(r+k-1)} x^r
				- \fr 1/{k(k-1)}\right]
\\
	&= n - 1 - \fr 1/2 \log n + \sum_{k=1}^\infty \fr B_{2k}/{2k(2k-1)}\left(1-\fr 1/{n^{2k-1}}\right),
		\qquad(n=1,2,3,\ldots).
\end{align*}

Our last step uses the following elementary summation formula for a recurring series,
\begin{align*}
	\sum_{r=0}^\infty  \fr \binom{r+k}{k}/{(r+k)(r+k-1)} x^r &=
				\fr 1/{k!} \sum_{r=0}^\infty  \fr {(r+k-2)!}/{r!} x^r
				=  \fr 1/{k!} D^{k-2}\sum_{r=0}^\infty  x^{r+k-2}
\\
	&= \fr 1/{k!} D^{k-2} \fr {x^{k-2}}/{1-x}
	= \fr 1/{k!} D^{k-2} \left(\fr {x^{k-2}-1}/{1-x}+\fr 1/{1-x}\right)
\\
	&= \fr {(k-2)!}/{k!} \fr 1/{(1-x)^{k-1}} = \fr {n^{k-1}}/{k(k-1)},   \qquad(k\ge2).
\end{align*}

After regrouping terms we finally obtain with De Moivre, for $n=1,2,3,\ldots$,
\begin{align} \nonumber
\log n! &= n\log n - \left\{ n - 1 - \fr 1/2 \log n + 
			\sum_{k=1}^\infty \fr B_{2k}/{2k(2k-1)}\left(1-\fr 1/{n^k} \right) \right\}
\\ \label{Deq}
	&= \left(n+\fr 1/2\right)\log n - n +
			  \left[1 - \sum_{k=1}^\infty \fr B_{2k}/{2k(2k-1)}\right]
			+ \sum_{k=1}^\infty \fr B_{2k}/{2k(2k-1)}\fr 1/{n^{2k-1}}.
\end{align}

\section{Quadrature of the circle by Wallis}

In order to determine the value of the constant expression $1-\sum_{k=1}^\infty \fr B_{2k}/{2k(2k-1)}$
in the previous equation (\ref{Deq}), we will use an infinite product for $\pi$ discovered via an 
ingenious induction by John Wallis in 1656, before the invention of the infinitesimal calculus.

	Taking limits in our now familiar Riemann sums, Wallis first computes $\int_0^1 x^r\,dx$ for a few
rational values of $r=\tfr p/q$, where $p$ and $q$ are small positive integers, then gives the 
general value of the integral as $\tfr 1/{(r+1)}$. From this, he can deduce with the integral binomial theorem
a few values of $\int_0^1 (x^p\pm x^q)^r\,dx$ where $p,q$ are rational and $r$ is a small positive integer,
and then he generalizes.

	A good number of known geometrical applications are given for parabolas, hyperbolas, spirals and
solids of revolution, but the real target is the ``quadrature of the circle", or the evaluation of
\begin{equation} \label{Wcircle}
	\fr \pi/4 = \int_0^1 \left(1-x^2\right)^{\fr 1/2} \, dx.
\end{equation}
Without the general binomial theorem (found later by Newton after reading this \cite{Nunn:1910}),
Wallis is forced to use interpolation between other integrals already computed.

\def\sq{\Box}

In the symmetric table of the reciprocals of $\int_0^1 (1-x^{\tfr 1/q})^p\,dx$ where $p,q=0,1,2,3,4,5$, 
he detects the figurate numbers \cite{ConwayGuy:1996}
whose values $\binom{p+q}{q}$ had been found earlier by Harriot and Fermat. What he is looking for corresponds to 
$p=q=\tfr 1/2$, so he inserts ``odd" lines and columns in the table for $p,q=-\tfr 1/2,\tfr 1/2,\tfr 3/2,\ldots$.
Their missing entries are then found by applying multiplication rules relating the values of
figurate numbers in the other ``even" lines and columns: Wallis knows intuitively that some general rules
are valid for all lines and columns which indeed contain $\fr \Gamma(p+q+1)/{(\Gamma(p+1)\Gamma(q+1))}$ 
\cite[p. 298]{Binet:1839}. The expanded table can thus be completely filled \cite{Nunn:1910,Osler:2010}, 
leaving a single unknown factor in the entries of the odd lines and odd columns, 
the ratio of the diameter squared to the circle area which is denoted by the small square ``$\sq$".

	Next, Wallis looks at the entries that he had found in the line corresponding for $q=\tfr 1/2$ to the 
successive multiples of one half, $p=-\tfr 1/2,\tfr 0/2,\tfr 1/2,\tfr 2/2,\tfr 3/2,\tfr 4/2,\tfr 5/2,\tfr 6/2,\ldots$~:
$$
	\fr 1/2 \sq,\quad  1,\quad \sq,\quad \fr 3/2,\quad  \fr 4/3 \sq,\quad  \fr 3\times 5/{2\times4},\quad
		  \fr {4\times 6}/{3\times5} \sq,\quad \fr {3\times 5\times7}/{2\times4\times6},\quad \cdots
 $$
The ratios of an odd/even column entry to the previous odd/even column entry obviously decrease as 
$(\tfr 2/1,\tfr 3/2,\tfr 4/3,\tfr 5/4,\ldots)$, and so should the ratios of one column to the previous one, giving
\begin{align*}
			\fr 1/{\tfr \sq/2} > \fr \sq/1 > \fr {\tfr 3/2}/{\sq} 
				&\quad\implies\quad  \sqrt{1+\fr 1/1} > \sq > \sqrt{1+\fr 1/2},
\\
			\fr {\tfr 3/2}/\sq > \fr {\tfr 4\sq/3}/{\tfr 3/2} > \fr {\tfr (3\times5)/{(2\times4)}}/{\tfr 4\sq/3} 
			&\quad\implies\quad
			\fr {3\cdot3}/{2\cdot4} \sqrt{1+\fr 1/3} > \sq > \fr {3\cdot3}/{2\cdot4} \sqrt{1+\fr 1/4}.
\end{align*}	% De_moivre 1730 Miscellanea Analytica p. 173
The other columns give similar inequalities, thus completing the quadrature of the circle (\ref{Wcircle})~:
$$
	 \sq := \fr 4/\pi = \fr 	{3\cdot3\cdot5\cdot5\cdots(2n-1)(2n-1)}/
						{2\cdot4\cdot4\cdot6\cdots(2n-2)(2n)\phantom{\,-0}} \theta_n,
			\quad \left(\sqrt{1+\fr 1/{2n-1}}>\theta_n>\sqrt{1+\fr 1/{2n}},n\ge1\right).
 $$

	Wallis shows in passing how to get the integer power sum formulas from the figurate numbers
values via Vandermonde's identity \cite[p. 140]{StedallWallis:2004}, but fails to notice in them the pattern
of recurring rational numbers that James Bernoulli would find fifty years later, using the same method.

Sharp criticisms from some scientists of Wallis' time were directed at his inductive method
based on informal analogies with arithmetic progressions, explained thus 
in Proposition 1 \cite[p. 13]{StedallWallis:2004}~:
\begin{quote}
``The simplest method of investigation, in this and various problems that follow,
is to exhibit the thing to a certain extent and to observe the ratios produced
and to compare them to each other; so that at length a general proposition may
become known by induction."
\end{quote}
He responded that this was appropriate during discovery, and we adopted his concept of limit (``the space will be less than any assigned quantity, and may therefore be taken as nothing") \cite[p. 145]{StedallWallis:2004}.

	Now the finite partial product can be expressed by factorials \cite{Cesaro:1880}
$$
     \fr \pi/2 \, \theta_n = \fr 1/{2n} \left[\fr2/1\times\fr4/3\times\fr6/5\times\ldots\times\fr2n/{2n-1}\right]^2
					= \fr 2^{4n}/{2n} \left[\fr (n!)^2/{(2n)!}\right]^2
 $$
and we obtain by taking logarithms
\begin{equation} \label{Stirling42}
	4 \log n! - 2\log (2n)! = \log 2n - 4n\log 2 + \log(\fr \pi/2\theta_n).
\end{equation}
	We also rewrite De Moivre's formula (\ref{Dform}) by introducing a constant C and a correction term $\delta_n$~:
\begin{align} \label{Stirling1}
	\log n! &= (n+\fr 1/2)\log n - n + C + \delta_n,
\\ \label{Stirling2}
	\log (2n)! &= (2n+\fr 1/2)\log 2n - 2n + C + \delta_{2n}.
\end{align}

	Comparing (\ref{Stirling42}) to the corresponding $(4,-2)$ linear combination of the last two equations 
(\ref{Stirling1}) and (\ref{Stirling2}) finally yields
$$
	C = \lim_{n\to\infty} \left(\log \sqrt{ 2\pi\theta_n } + \delta_{2n} - 2 \delta_n\right) = \log \sqrt{ 2\pi},
 $$
assuming, as may be easily shown \cite{Cesaro:1880,Rouche:1890}, that the correction term 
$\delta_n$ vanishes as $n\to\infty$.

\bigskip

\section{Stirling's derivation of the original series}

While the proof by de Moivre could be purely analytical since he knew the formula, the original
proof by Stirling looks more experimental, like a sophisticated version for series \cite{Stirling:1749} of the
``inductive" method of pattern discovery used earlier by John Wallis for the quadrature of the circle.

	The problem considered in Proposition XXVIII of ``Methodus Differentialis" is the evaluation of
a finite sum of logarithms in arithmetic progression of step $2h$,
$$
	S = \log(x+h) + \log(x+3h) + \ldots + \log(z-h),
 $$
which Stirling expresses as a telescoping sum of differences of a ``finite integral" $F(z)$
(see \cite{Prouhet:1864}),
$$
	S = \left[F(x+2h)-F(x)\right] + \left[F(x+4h)-F(x+2h)\right] + \ldots + \left[F(z)-F(z-2h)\right].
 $$
Obviously, it is sufficient to choose the function so that $F(z)-F(z-2h)=\log(z-h)$, and Stirling simply
indicates that this is satisfied by
$$
	F(z) = \fr {z\log z-z}/{2h} - \fr 1/{12} \fr h/{z} + \fr 7/{360}\fr h^3/{z^3}
	- \fr 31/{1260}\fr h^5/{z^5} + \fr 127/{1680}\fr h^7/{z^7} - \fr 511/{1188}\fr h^9/{z^9} + \ldots.
 $$
The constants in the series must be computed from linear equations formed by alternate binomial coefficients
of odd order~:
\begin{align*}
	- \fr 1/{3\times4} &= a_1,
\\
	- \fr 1/{5\times8} &= a_1 + 3a_3,
\\
	- \fr 1/{7\times12} &= a_1 + 10a_3 + 5a_5,
\\
	- \fr 1/{9\times16} &= a_1 + 21a_3 + 35a_5 + 7a_7,
\\
	- \fr 1/{11\times20} &= a_1 + 36a_3 + 126a_5 + 84a_7 + 9a_9.
\end{align*}
This ends the short proof, to which a hint to a connection with an ``area correction" is added.

	Next, in the very important example II, Stirling indicates that the logarithm of $n!$ can be 
approximated by computing three or four terms of the series in $F(n+\tfr1/2)$, then adding half the
logarithm of the circumference of a unit circle, equal to $0.39908 99341 79$ in base $10$.
This constant $\log \sqrt{2\pi}$ comes from the Wallis product that he used earlier to evaluate 
$\binom{2n}{n}$ in Proposition XXIII.

		Indeed, we can verify the previous equations as follows. First, since the sum of logarithms is
a Riemann sum, a leading term $\int \log z dz=z\log z - z$ is expected \cite{Prouhet:1864}, 
and a correcting series with undetermined coefficients $a_k$ can be added, giving as candidate
for the telescoping finite integral
$$
	F(z) = \fr {z\log z - z}/{2h} + \sum_{k=0}^\infty a_k \fr h^k/{z^k}.
 $$
Secondly, we set $z=x+h$ and $t=\tfr h/x$, which eliminates $x$ in
\begin{align*}
	&F(z) - F(z-2h) - \log(z-h)
\\
	&= F(x+h) - F(x-h) - \log x
\\
	&=  \quad (x+h) \fr { \log x + \log(1 + \tfr h/x) - 1 }/{2h} + \sum_{k=0}^\infty a_k \fr h^k/{(x+h)^k}
\\
	& \quad - (x-h) \fr { \log x + \log(1 - \tfr h/x) - 1 }/{2h} - \sum_{k=0}^\infty a_k \fr h^k/{(x-h)^k}
	  	  - \log x
\\
	&= \fr 1/{2t} \log \fr {1+t}/{1-t} + \fr 1/2\log(1-t^2) - 1
		+ \sum_{k=1}^\infty a_k t^k\left[\fr 1/{(1+t)^k}-\fr 1/{(1-t)^k}\right]
\\
	&= \sum_{j=1}^\infty \fr {1-(-1)^j}/2 \fr t^{j-1}/j
		 - \sum_{j=1}^\infty \fr t^{2j}/{2j}
		+ \sum_{k=1}^\infty a_k t^k \fr D^{k-1}/{(k-1)!} \left[\fr (-1)^{k-1}/{1+t}-\fr 1/{1-t}\right]
\\
	&= \sum_{j=1}^\infty \left[ \fr t^{2j}/{2j+1} - \fr t^{2j}/{2j} \right]
		+ \sum_{k=1}^\infty a_k \fr t^k/{(k-1)!} 
			\sum_{j=k-1}^\infty  \fr {j!\; t^{j-k+1}}/{(j-k+1)!} \left[(-1)^{j-k+1}-1\right]
\\
	&= - \sum_{j=1}^\infty \fr t^{2j}/{2j(2j+1)}
		+ \sum_{j=0}^\infty t^{j+1}\sum_{k=0}^j \binom{j}{k} a_{k+1} \left[ (-1)^{j-k} - 1 \right].
\end{align*}
The right hand side will vanish identically if (for $j=2n$ in the second sum)
$$
	0 = \sum_{k=0}^{2n} \binom{2n}{k} a_{k+1} \left[ (-1)^{2n-k} - 1 \right]
	= - 2 \, \sum_{k=0}^{n-1}\binom{2n}{2k+1} a_{2k+2};
 $$
and also if (for $j=n$ in the first sum and $j=2n-1$ in the second sum)
$$
	-\fr 1/{2n(2n+1)} = - \sum_{k=0}^{2n-1} \binom{2n-1}{k} a_{k+1} \left[ (-1)^{2n-1-k} - 1 \right]
	 = 2 \, \sum_{k=0}^{n-1} \binom{2n-1}{2k} a_{2k+1}, \quad(n=1,2,\ldots).
 $$
Clearly, as indicated by Stirling, the $a_{2k}$ must vanish for $k>0$ while the $a_{2k+1}$ must satisfy
$$
	- \fr 1/{(2n+1)(4n)} = \sum_{k=0}^{n-1} \binom{2n-1}{2k} a_{2k+1}, \qquad(n=1,2,\ldots).
 $$
The last recurrence relation is equivalent to 
$$
	 \fr 1/2 + \sum_{k=1}^{n} \binom{2n+1}{2k} 2k(2k-1)a_{2k-1} =0, \qquad(n=1,2,\ldots).
 $$
Now, as stated by James Bernoulli, with $B_0=1$, $B_1=-\tfr 1/2$, and $B_{2k+1}=0$ for $k>0$,
$$
	1 - \fr 1/2 + \sum_{k=1}^{n} \binom{2n+1}{2k} B_{2k} 
	= \sum_{k=0}^{2n} \binom{2n+1}{k} B_{k} =0, \qquad(n=1,2,\ldots).
 $$
Thus, as seen by De Moivre while reading the letter of Stirling,
$$
	a_{2k-1} = \fr {B_{2k}}/{2k(2k-1)}, \qquad(k=1,2,\ldots).
 $$
The equivalence of De Moivre's and Stirling's versions of the series can now be shown easily.

\section{Semi-convergence}
De Moivre had only the six $B_{2k}$ numbers that we listed above, from his copy of Ars Conjectandi.
Although he reformulated the recurrence rule given by James Bernoulli to determine more $B_{2k}$, he never computed 
enough terms to notice the divergence of the correcting series for $\log n!$. 
This was proven later in 1763 by Thomas Bayes in a posthumous letter published by the Royal Society in London~:

\begin{quote}
``It has been asserted by some eminent mathematicians, that the sum of the logarithms of the numbers
$1\cdot2\cdot3\cdot4\cdot$\&c. to $n$, is equal to 
$$
	\fr 1/2\log C + \left(n+\fr 1/2\right)\log n
 $$
lessened by the series 
$$
	n - \fr 1/{12n} + \fr 1/{360n^3} - \fr 1/{1260n^5} + \fr 1/{1680n^7} - \fr 1/{1188n^9} + \&c.
 $$
if $C$ demote the circumference of a circle whose radius is unity. (.........)

Whereforeatlength the subsequent terms of this series are greater than the preceding ones,
and increase in infinitum, and therefore the whole series can have no ultimate value whatsoever.
Much less can that series have any ultimate value, which is deduced from it by taking $n=1$,
and is supposed to be equal to the logarithm of the square root of the periphery of a circle
whose radius is unity\footnote{In 1755, Euler had stated the divergence in
		 $1-\fr B_2/{1\cdot2}-\fr B_4/{3\cdot4}-\cdots = \log\sqrt{2\pi}$
		 \cite[C. VI, \S157, p. 465]{Euler:1755}.}."
\end{quote}
\noindent
The Reverend Bayes was careful in using the term ``whole series", since he had noted that the first terms
of the series gave a good approximation, in particular when $n$ is large. 

In fact the best number of terms to use has been determined to be near the integer part of 
$\pi n$, after which terms stop decreasing in absolute value \cite{Legendre:1811,Schaar:1848,Limbourg:1859}.
Moreover, the rest of the series is bounded in absolute value by the first 
neglected term and has the same sign, which alternates from term to term following the sign of the 
Bernoulli numbers $B_{2k}$ \cite{Cauchy:1841}. This is an ``enveloping series" whose successive partial sums
encompass a specific value which can be assigned in a certain sense to the whole series.
For example, George Boole writes in 1860 \cite[p. 87]{Boole:1860}:
\begin{quote}
	``Hence for large values of $n$ we may assume
$$
	1\cdot2\cdot\cdots n=\sqrt{2\pi n} \left(\fr n/e\right)^n,
 $$
the ratio of the two members tending to unity as $n$ tends to infinity. And speaking generally it
is with the ratios, not the actual values of functions of large numbers, that we are concerned."
\end{quote}

	In order to obtain a formula valid according to our modern standards, it is necessary to
replace the divergent infinite series by its partial sum, and find an expression for the 
neglected rest of the series. This is usually done first for the Euler-Maclaurin formula,
and then applied generally to $\log \Gamma(z)$, where $z$ is a complex number outside the 
negative real axis \cite{Henrici:1977}. But various shortcuts have been found when $z$ is 
restricted to the positive real axis, using properties of the $\Gamma$ function and infinite series
\cite{Glaisher:1878,GlaisherMansion:1879,Cesaro:1880,Texeira:1891,Godefroy:1903,Nielsen:1906}.
In \cite{Schaar:1848}, Schaar replaced the series involving the Bernoulli numbers by an integral,
using the Euler partial fraction decomposition of the cotangent function,
$$
	\fr 1/{e^x-1} = \fr 1/x - \fr 1/2 + \sum_{k=1}^\infty \fr 2x/{x^2+4k^2\pi^2},
 $$
deduced from the partial fractions of $\tfr 1/{(t^{2n}-1)}$ where $t=1+\tfr x/{2n}$
and $e^x=\lim (1+\tfr x/n)^n$ \cite{Bernoulli:1685}.
From clever manipulations of the integral for $\Gamma(a+1)$, he finally obtains, for $a>0$ and $m\ge0$,
\begin{align*}
	\log \Gamma(a+1) &= \log\sqrt{2\pi a} + a(\log a - 1)
				+ \sum_{k=1}^m \fr B_{2k}/{2k(2k-1)} \fr 1/{a^{2k-1}}
\\
	&\qquad\qquad	- (-1)^m \fr 1/\pi \int_0^\infty \fr x^{2m}/{1+x^2} \log(1-e^{-2\pi ax})\,dx.
\end{align*}

\section{Errors in the table of logarithms of $n!$ by Abraham de Moivre}

The 1756 edition of the first De Moivre book \cite{DeMoivre:1718} contains in Appendix IV on page 333
``A correct Table of the Sums of Logarithms, from the Author's Supplement to his Miscellanea Analytica"
that we compare to the original table of $\log (10n)!$ for $1\le n\le90$
on pages 103--104 in the 1730 ``Miscellanea Analytica". 
Spaces between digits in the following table section highlight differences which are also represented
in the last column.
\def\z{\phantom{0}}
\def\dfive{  10^{-5}     }
\def\dseven{  6\cdot10^{-7}  }
\def\dnine{ 9\cdot10^{-9} }
\def\dten{ 27\cdot10^{-10}}

$$
\begin{array}{|r|l|l|l|}
\hline
	n	& \quad \log(n!)\quad(1730)	&\quad \log(n!)\quad(1756)	&  \quad(1730) - (1756)
\\
\hline
 10 &\z\z6.55977\z 303287678      &  \z\z6.55976\z 303287678\z\z  & +\dfive
\\
 20 & \z18.38613\z 461687770      &  \z18.38612\z 461687770\z\z   & +\dfive
\\
 30 & \z32.42367\z 007492572      &  \z32.42366\z 007492572\z\z   & +\dfive
\\
 40 & \z47.91165\z 506815591\z\z  &  \z47.91164\z 506815591\z\z   & +\dfive
\\
 50 & \z64.48308\z 487247209\z\z  &  \z64.48307\z 487247209\z\z   & +\dfive
\\
 60 & \z81.92018\z 484939024\z\z  &  \z81.92017\z 484939024\z\z   & +\dfive
\\
 70 &  100.07841\z 503568004\z\z  &   100.07840\z 503568004\z\z   & +\dfive
\\
 80 &  118.85473\z 71\z 2249966\z &   118.85472\z 77\z 2249966\z  & +\dfive\,-\,\dseven
\\
 90 &  138.17194\z 51\z 9001086\z &   138.17193\z 57\z 9001086\z  & +\dfive\,-\,\dseven
\\
100 &  157.97001\z 30\z 5471585\z &   157.97000\z 36\z 5471585\z  & +\dfive\,-\,\dseven
\\
110 &  178.20092\z 70\z 4487008\z &   178.20091\z 76\z 4487008\z  & +\dfive\,-\,\dseven
\\
120 &  198.82540\z 32\z 4721977\z &   198.82539\z 38\z 4721977\z  & +\dfive\,-\,\dseven
\\
130 &  219.81070\z 25\z 5614815\z &   219.81069\z 31\z 5614815\z  & +\dfive\,-\,\dseven
\\
140 &  241.12911\z 93\z 8869689\z &   241.12910\z 99\z 8869689\z  & +\dfive\,-\,\dseven
\\
150 &  262.75690\z 28\z 1092616\z &   262.75689\z 34\z 1092616\z  & +\dfive\,-\,\dseven
\\
160 &  284.67346\z 56\z 4068298\z &   284.67345\z 62\z 4068298\z  & +\dfive\,-\,\dseven
\\
170 &  306.66079\z 13\z 9482847\z &   306.66078\z 19\z 9482847\z  & +\dfive\,-\,\dseven
\\
180 &  329.30298\z 08\z 238\z 9393 &  329.30297\z 14\z 247\z 9393 & +\dfive\,-\,\dseven\,-\,\dnine
\\
190 &  351.98589\z 92\z 366\z 3535 &  351.98588\z 98\z 339\z 3535 & +\dfive\,-\,\dseven\,+\,\dten
\\
200 &  374.89689\z 80\z 427\z 4044 &  374.89688\z 86\z 400\z 4044 & +\dfive\,-\,\dseven\,+\,\dten
\\
\hline
\end{array}
 $$
Clearly, the systematic errors which recur from one line to the next are not typographic errors.
One can see that there was in 1730 an error of $+1$ unit in the fifth decimal of $\log(10!)$, which affected
all other entries. And for $\log(80!)$, "1" replaced a "7" in the seventh decimal, also affecting
all subsequent entries below. Then a most probable inversion of the ninth and tenth decimals 
("74" instead of "47") in $\log(180!)$ was also propagated to the rest of the table.

Based on this data and its title, one can speculate that De Moivre simply added the logarithms of 
natural numbers in sequence to compute the logarithms of the factorials, without using as a check
an approximating formula, another order of additions, or sums of digits modulo $9$ 
(``casting out the nines" \cite{Luroth:1900,Chrystal:1904,ConwayGuy:1996}, which does not detect 
digit permutations).

\section{Stirling's series in ``The Doctrine of Chances" by Abraham de Moivre, 1756}
The 1756 edition of the first De Moivre book \cite{DeMoivre:1718} also contains
on page 334 the following presentation of the original Stirling series after the
corrected table of logarithms.

\begin{quote}

``If we would examine these numbers, or continue the Table farther on, we have that excellent
rule communicated to the Author by Mr. James Stirling; published in his
Supplement to the ``Miscellanea Analytica" and by Mr. Stirling himself in his
``Methodis Differentialis", Prop. XXVIII.

Let $z-\tfr 1/2$ be the last term of any Series of the Natural Numbers 
$1,2,3,4,\ldots,z-\tfr 1/2$; $a=.43428448190325$ the reciprocal of Neper's Logarithm of 10:
\\
Then three or four terms of this series 
$$
	z\log z - az - \fr a/{2\cdot12\cdot z} + \fr 7a/{8\cdot360z^3}
	- \fr 31a/{32\cdot1260z^5} + \fr 127a/{128\cdot1680z^7} - \&c.
 $$
added to $0.399089934179$ which is half the Logarithm of a Circumference whose Radius is Unity,
will be the Sum of the Logarithms of the given Series; or the Logarithm of the Product
$1\times2\times3\times4\times5\times\cdots\times(z-\tfr 1/2)$.

The coefficients of all the terms after the first two being formed as follows. Put
\begin{align*}
	- \fr 1/{3\cdot4} &= A
\\
	- \fr 1/{5\cdot8} &= A + 3 B
\\
	- \fr 1/{7\cdot12} &= A + 10 B + 5 C
\\
	- \fr 1/{9\cdot15} &= A + 21 B + 35 C + 7 D
\\
	- \fr 1/{11\cdot20} &= A + 36 B + 126 C + 84 D + 9 E.
\end{align*}

\medskip

In which the Numbers $1,\,1,\,1,\,$\&c., $\quad 3,\,10,\,21,\,36,\,$\&c., $\quad 5,\,35,\,126,\,$\&c.
that multiply $A,B,C,\&c.$ are the alternate Unciae of the odd Powers of a Binomial.

Then the Coefficients of the several Terms will be
$(\tfr 1/2)\times A=-\tfr 1/{(2\cdot12)}$,\\
$(\tfr 1/2)^3\times B=-\tfr 7/{(8\cdot360)}$,
$(\tfr 1/2)^5\times C=-\tfr 31/{(32\cdot1260)}$, \&c. 

See the general Theorem and Demonstration in Mr. Stirling's Proposition quoted above."

\end{quote}

\bigskip

\section{Conclusion}

The errors in the table of logarithms were in fact a godsend, since they prompted the 
letter of Stirling; and by chance De Moivre then worked not only on correcting his additions but 
also on finding a new proof of the original formula (\ref{Sform}).
His simpler version (\ref{Dform}) of the series is now universally adopted because so much more
properties are known for the Bernoulli numbers than for the coefficients of Stirling, after the
appearance in 1755 of the  revolutionary ``Institutiones Calculi Differentialis" by Euler
\cite{Euler:1755} -- another lucky break for De Moivre.

Stirling himself admitted later gracefully in a September 1730 letter to Cramer that De Moivre's version 
was simpler than his own \cite{Dutka:1991}.
He also suggested in 1738 the formula (which must then be De Moivre's version) to Euler as an example 
for his new summation method using differential quotients and Bernoulli numbers,  while making him aware
of the existence of the similar Maclaurin's summation method \cite{Fowler:2000}.

Now De Moivre first discovered, proved geometrically, and published the trigonometric formula 
$(\cos\theta+i\sin\theta)^n=\cos n\theta+i\sin n\theta$.
While Euler's reformulation and proof via his own relation $e^{i\theta}=\cos\theta+i\sin\theta$ is now
universally known and preferred, nobody (certainly not Euler!)
would suggest switching priority from De Moivre to Euler for that new proof and new notation.

Likewise indeed, James Stirling must keep his indisputable priority for the discovery, proof, 
and publication of Stirling's series (\ref{Sform}), including of course the explicit constant
$\log\sqrt{2\pi}$.

%=================================================================================================
\bibliographystyle{plain}

%=================================================================================================

\end{document}